\documentclass[a4paper,10pt]{article}
\usepackage[left=44mm,right=40mm,top=60mm,bottom=52mm]{geometry}
\usepackage{amsmath}
\usepackage{amssymb}

\newtheorem{theorem}{Theorem}

\newtheorem{lemma}{Lemma}
\newtheorem{proposition}{Proposition}
\newtheorem{remark}{Remark}

\begin{document}

\section*{On maximum of Gaussian process \\
with unique maximum point of its variance\footnote{Partially
supported by Russian Science Foundation, grant 14-49-00079 and
SNSF grant 200021-175752} }

\begin{flushright}
\textbf{E. HASHORVA}\\
\textit{University of Lausanne\\
Enkelejd.Hashorva@unil.ch}\\\textbf{S. G. KOBELKOV}\\
\textit{Lomonosov Moscow State University\\
sergeyko81@gmail.com }\\
\textbf{V. I. PITERBARG}\\
\textit{Lomonosov Moscow State University,\\
National Research University "MPEI",\\
Scientific Research Institute of System Development "NIISI RAS"\\
piter@mech.math.msu.su}
\end{flushright}

UDC 519.218\newline

{\small \textbf{Key words:} non-stationary process, Gaussian
process, maximum, Pickands' method, double sum method. }

\begin{center}
{\small \textbf{Abstract} }
\end{center}

{\small Gaussian random processes which variances reach theirs maximum
values at unique points are considered. Exact asymptotic behaviors of
probabilities of large  absolute maximums of theirs trajectories have been
evaluated using Double Sum Method under  the widest possible conditions. }

\section{ Introduction. Preliminaries.}

This note is a generalization of \cite{PP}. Our aim is to show the maximum
capability of the Pickands' Double Sum Method for asymptotic behavior of the
maximum tail distribution for Gaussian stationary process, see \cite%
{pickands}, with corrections in \cite{pit1}. This method has been
generalized to Gaussian random fields, \cite{book}, where stationary fields
with power like behavior of the correlation function at zero are considered
as well as fields with a similar behavior of the correlation function at the
unique maximum point of variance and with power like behavior of it near the
point. However, while the power behavior of the correlation function, with
possible light generalization to regular variation of it, \cite{pit1}, is
quite essential for the Pickand's method, the required in \cite{PP}, \cite%
{book} power behavior of the variance looks somewhat artificial. In the
present note we give the widest possible conditions on the variance and on
the correlation function under which the Double Sum Method still works. Note
also that in the recent article \cite{HDL} it is proved that in the
non-stationary case the variance behavior does not need to be power but can
just be regularly varying.

Let $X(t),$ $t\in \lbrack -S,S],$ be a zero mean a.s. continuous Gaussian
process with covariance function $r(s,t),$ denote $\sigma ^{2}(t)=r(t,t)$.
Here we study the asymptotic behavior of the probability
\begin{equation}
P([-S,S];u):=P(\max_{t\in \lbrack -S,S]}X(t)>u)  \label{P(S,u)}
\end{equation}%
as $u\rightarrow \infty .$ We assume that $\sigma (t)$ reaches its absolute
maximal value only at zero, since in the case of another point of the
absolute maximum one can simply shift the time.

Assume the following.

\begin{description}
\item[ $\mathbf{A1:}$] Suppose that $X$ has a.s.\ continuous sample paths.
\end{description}

In particular, the above assumption is satisfied under the following
standard H\"{o}lder condition, namely for some positive $\Gamma$ and $\gamma$%
,
\begin{equation}
E (X(t)-X(s))^{2}\leq \Gamma|t-s|^{\gamma},\ \ s,t\in\lbrack-S,S].
\label{A1}
\end{equation}
Under this condition there exits an a.s. continuous version of $X$. Here, in
contrast of \cite{PP}, see also \cite{lectures}, \cite{book}, we do not
assume (\ref{A1}).

\begin{description}
\item[$\mathbf{A2:}$] $\sigma(t)$ reaches its global maximum on $[-S,S]$
only at  $0\mathbf{\ }$and $\sigma(0)=1.$ Moreover, there exist finite or
infinite limits%
\begin{equation}
\lim_{t\downarrow0,s\downarrow0}\frac{1-\sigma^{2}(t)}{1-r(s,t)}\in
\lbrack0,\infty]\ \ \text{and }\lim_{t\uparrow0,s\uparrow0}\frac{1-\sigma
^{2}(t)}{1-r(s,t)}\in\lbrack0,\infty].  \label{A2}
\end{equation}
\end{description}

Note again that the above specification of the location and the maximal
value of $\sigma (t)$ is just for convenience. Note further that $\mathbf{A2}
$ implies $r(s,t)\leq 1,\forall s,t\in \lbrack -S,S]$ with equality holding
only for $s=t=0.$

Denote by $\rho(s,t):=r(s,t)/\sigma(s)\sigma(t)$, the correlation function
of $X$.

\begin{description}
\item[$\mathbf{A3:}$ (Local stationarity at $0$)] There exists a covariance
function $\rho (t)$ of a stationary process such that
\begin{equation*}
\lim_{s,t\rightarrow 0}\frac{1-\rho (s,t)}{1-\rho (t-s)}=1.
\end{equation*}

\item[$\mathbf{A4:}$] For $\rho$ from $\mathbf{A3}$, there exist a positive
function $q(u)$ and a function $h(t),$ $h(t)>0$ for all $t\neq0$ such that
\begin{equation}
\lim_{u\rightarrow\infty}u^{2}(1-\rho(q(u)t\mathbf{))=}h(t)  \label{A3}
\end{equation}
uniformly over $t\in\lbrack-\varepsilon,\varepsilon]$ for some $%
\varepsilon>0.$
\end{description}

Notice that $\rho (0)=1$ and $\mathbf{A3}$ imply that $\rho $ is continuous,
hence $q(u)\rightarrow 0$ as $u\rightarrow \infty ,$ and therefore (\ref{A3}%
) is fulfilled uniformly over any compact set. Furthermore, it also follows
from (\ref{A3}) that for any positive $s,t,$
\begin{equation}
\lim_{u\rightarrow \infty }\frac{1-\rho (q(u)t\mathbf{)}}{1-\rho (q(u)s%
\mathbf{)}}=\frac{h(t)}{h(s)},  \label{e2}
\end{equation}%
which implies, by definition, the regular variation at zero of $1-\rho (t),$
\cite{BGT}. The index of the regular variation, say $\alpha ,$ is positive,
and $h(t)=t^{\alpha }.$ Indeed, if $\alpha <0,$ $\rho (t)$ is not continuous
at zero, if $\alpha =0,$ $h(t)=1$ for all $t>0$ and $h(t)=0$ for $t=0,$ so
it is not continuous again. Further, if $\alpha >2$ it follows from $\mathbf{%
A3}$ and $\mathbf{A4}$ that $\rho ^{\prime \prime }(t)\equiv 0$ which
contradicts the positive definiteness of $\rho $. Consequently, we have that
$\alpha \in (0,2].$ As well, the same is valid for $\alpha =2$ and $%
t^{-2}(1-\rho (t))\rightarrow 0$. Thus, assumption $\mathbf{A4}$ is
equivalent to the corresponding assumption in \cite{pit1}, and therefore
this condition is crucial for our method, the Double Sum Method. Thus we
have,
\begin{equation}
1-\rho (t)\ \text{\emph{is regularly varying at zero with index} }\alpha \in
(0,2].  \label{alpha}
\end{equation}%
Further, since $1-\rho (t)=\ell (t)t^{\alpha },$ where $\ell (t)$ is slowly
varying function at zero, we have,%
\begin{equation*}
q(u)=(1-\rho )^{\leftarrow }(u^{-2}),
\end{equation*}%
where "$^{\leftarrow }$" means the generalized inverse. Now using Theorems
1.5.12, 1.5.13 (de Bruijn Lemma), and Proposition 1.5.15, \cite{BGT}, we get
that
\begin{equation}
q(u)\sim u^{-2/\alpha }\ell ^{\#}(u^{-2})^{1/\alpha }  \label{q(u)}
\end{equation}%
as $u\rightarrow \infty .$ In our notation $\sim $ stands for
asymptotic equivalence, and $\ell ^{\#}$ is the de Bruijn
conjugate of $\ell .$ In view of (\ref{e2}), we have that
(\ref{A3}) holds for any $q^{\prime }$ such that
$\lim_{u\rightarrow \infty }q(u)/q^{\prime }(u)=1.$ Consequently,
since $q$ is regularly varying at infinity, without loss of
generality we assume hereinafter that $q$ is monotone. \newline
Note that the slowly varying function $\ell ^{\#}$ can be often
explicitly
calculated, see Bojanic and Seneta Theorem 2.3.3 and Corollary 2.3.4, \cite%
{BGT}. For example, if
\begin{equation*}
\frac{\ell \left( u^{-2}\ell (u^{-2})\right) }{\ell (u^{-2})}\rightarrow 1\
\text{as }u\rightarrow \infty \text{,}
\end{equation*}%
then $\ell ^{\#}\sim 1/\ell .$

In Section 2 we repeat the results from \cite{pit1} in this new conditions.
In Section 3 the main result of the paper is presented. In short Section 4
we present two examples to demonstrate the generality of our result.

\section{Stationary processes}

In this section we assume that $X(t),$ $t\in \lbrack 0,S],$ is a stationary
Gaussian process with mean zero, unit variance and covariance function $\rho
$ described above. We formulate here for convenience the results from \cite%
{pit1} with some obvious further generalizations.

\begin{lemma}
\label{lemma_pickands} If $\mathbf{A1}$ and $\mathbf{A4}$ hold, then for any
$T>0,$%
\begin{equation*}
P([0,q(u)T];u)=(1+\gamma (u))H_{\mathbf{\alpha }}(T)\Psi (u),
\end{equation*}%
with $\gamma (u)\rightarrow 0$ as $u\rightarrow \infty ,$ where $\alpha \in
(0,2]$ is defined in (\ref{alpha}), by the arguments below $\mathbf{A4,}$
\begin{equation*}
H_{\alpha }(T)=\mathbf{E}\exp (\max_{[0,T]}\chi (t)),\ \
\end{equation*}%
and $\chi (t\mathbf{)}$ is a Gaussian process with continuous trajectories, $%
\chi (0)=0,$ and
\begin{equation*}
\mathrm{var}(\chi (t)-\chi (s))=2h(|t-s|\mathbf{),\ E}\chi (t)=-h(t).
\end{equation*}
\end{lemma}

\begin{theorem}
\label{Pickands} Suppose that the conditions of Lemma 1 hold. Let
furthermore $\rho (t)<1$ for all $t>0$. Then for any $E\subset \mathbb{R},$
a bounded closure of an open set,%
\begin{equation*}
P(E;u)=\mathrm{mes}(E)\mathcal{H}_{{\alpha }}\frac{\Psi (u)}{q(u)}(1+o(1)),\
\ \ u\rightarrow \infty ,
\end{equation*}%
as $u\rightarrow \infty ,$ where
\begin{equation*}
\mathcal{H}_{\alpha }=\lim_{T\rightarrow \infty }T^{-1}H_{\alpha }(T)\in
(0,\infty ).
\end{equation*}%
This assertion holds even if $E=E(u)$, provided there exist segments $%
E^{-}(u),E^{+}(u)\subset \mathbb{R}$ such that $E^{-}(u)\subset E\subset
E^{+}(u)$ with $\lim_{u\rightarrow \infty }\mathrm{mes}(E^{-}(u))/q(u)=%
\infty ,$ and for some $\delta \in (0,1/2)$, we have $\mathrm{mes}%
(E^{+}(u))e^{-\delta u^{2}}\rightarrow 0$ as $u\rightarrow \infty .$
\end{theorem}

\section{Gaussian processes with the unique maximum point of variance}

In this section we consider a centered non-stationary Gaussian process $X(t),
$ $t\in \lbrack -S,S]$. In view of $\mathbf{A4},$ it follows from $\mathbf{A2%
}$ that there exists the limit
\begin{equation}
\lim_{u\rightarrow \infty }u^{2}(1-\sigma ^{2}(q(u)t)\mathbf{)=}h_{1}(t)\in
\lbrack 0,\infty ].  \label{cond2}
\end{equation}%
Notice that the limit relations in $\mathbf{A2}$ follow from (\ref{cond2})
as well. The limit $h_{1}(t)$ can be equal to zero, it can be positive and
finite, it can be equal to infinity. These assertions do not change for any
other $t$ of the same sign, that is, from the same half-line. We say that
\emph{the stationary-like case} takes place if the limit equals zero for all
$t$, see discussion below. If the limit is equal to infinity, we shall refer
to the \emph{Talagrand-like case}, in this case for any set $S$ containing
zero
\begin{equation*}
P(S;u)\sim {P}(X(0)>u),\ \ u\rightarrow \infty ,
\end{equation*}%
see the proof below. Talagrand has shown this for general Gaussian processes
and under the most general conditions, see \cite{book} for references and
discussions. Finally, for non-zero and non-infinity $h_{1}(t)$, the third
case is called \emph{the transition case}. Since we do not assume that $%
\sigma $ is symmetric with respect to zero, consideration of left and right
limits in (\ref{cond2}) may has a combination of three cases above. For
instance, $h_{1}(t)=\infty $ for any $t\in \lbrack -S,0)$ and $h_{1}(t)\in
(0,\infty )$ for any $t\in (0,S]$. In the latter case, the arguments given
for $1-\rho (q(u)t),$ imply that for some $\beta \geq 0,$ $1-\sigma
^{2}(t),\ t>0$ is regularly varying at $0$ and moreover $%
h_{1}(t)=h_{1}(1)t^{\beta }.$ Since
\begin{equation*}
\lim_{u\rightarrow \infty }\frac{u^{2}(1-\sigma ^{2}(q(u)t)\mathbf{)}}{%
u^{2}(1-\rho (q(u)t)\mathbf{)}}=\frac{h_{1}(t)}{h(t)},
\end{equation*}%
we conclude that $\alpha =\beta $ and further, the regularly varying
functions $1-\sigma ^{2}(t)$ and $1-\rho (t)$ have to be equivalent up to a
positive constant, namely we have
\begin{equation*}
\lim_{t\downarrow 0}\frac{1-\sigma ^{2}(t)}{1-\rho (t)}=\frac{h_{1}(1)}{h(1)}%
>0.
\end{equation*}

Now we formulate two general results for all described above types of
behavior of $\sigma(t)$. The first one is a standard local lemma of Double
Sum Method, a generalization of Lemma \ref{lemma_pickands}, see \cite{book},
\cite{lectures}.

\begin{lemma}
\label{PP-lemma} Under the assumptions $\mathbf{A1}$ -- $\mathbf{A4}$, for
any $T>0$,%
\begin{equation*}
P([0,q(u)T];u)=P_{\alpha}^{+}(T\mathbf{)}\Psi(u)(1+o(1))
\end{equation*}
and%
\begin{equation*}
P([-q(u),q(u)T];u)=P_{\alpha}(T\mathbf{)}\Psi(u)(1+o(1)),
\end{equation*}
as $u\rightarrow\infty,$ where%
\begin{equation*}
P_{\alpha}^{+}(T)=E\max_{t\in\lbrack0,T]}e^{\chi_{1}(t)},\quad
P_{\alpha}(T)=E\max_{t\in\lbrack-T,T]}e^{\chi_{1}(t)},
\end{equation*}
with $\chi_{1}(t)=\chi(t)-h_{1}(t)$ for $h_{1}(t)<\infty$, and $\chi_{1}(t)=0
$ for $h_{1}(t)=\infty$.
\end{lemma}

The proof of this lemma is a simple repetition of the proof of Lemma 6.1
\cite{book} by using the assumptions $\mathbf{A1-A4}$ and the relation (\ref%
{cond2}). The case $h_{1}(t)=\infty$ can be treated by similar arguments.
Note that in the Talagrand case, the detailed consideration of the weak
convergence in $C([-T,T])$ of the process
\begin{equation*}
\chi_{u}(t)=u(X(q(u)t)-u)+w
\end{equation*}
given $X(0)=u-w/u$, can be restricted to $C([-T,0])$ for $t>0$, or to $%
C([0,T])$ for $t<0$. It can be proved that in these cases given $X(0)=u-w/u,$
\begin{equation*}
\max_{t\in\lbrack-T,T]}\chi_{u}(t)\rightarrow\max_{t\in\lbrack-T,0]}(%
\chi(t)-h_{1}(t))\ \ \text{as }u\rightarrow\infty,
\end{equation*}
weakly for $t>0$, and similar convergence hold for $t<0.$ If $h_{1}(t)=\infty
$ for all non-zero $t,$ the above weak limit is equal to $0.$

The next result concerns the extraction of an informative parameter set
depending on the level $u$, which provides the required asymptotic behavior.
Consider the set

\begin{equation}
B_{u}=\left\{ t\mathbf{:}1-\sigma^{2}(t){\leq}u^{-2}\log^{A}u\right\} ,\ \
A>1.  \label{B_u}
\end{equation}

\begin{lemma}
\label{extracting} If $X$ is a centered Gaussian process satisfying $\mathbf{%
A1}-\mathbf{A4}$, then for any $E\subset [-S,S]$, which is a closure of a
bounded open set containing zero, and for any $B\in(1,A)$, we have
\begin{equation*}
P(E;u)= P(E\cap B_{u};u)\left( 1+O\left( e^{-\log^{B}u}\right) \right)
\end{equation*}
as $u\rightarrow\infty.$
\end{lemma}

\textbf{Proof: } By \textbf{A1} $X$ has bounded sample paths almost surely.
Then the Borell-TIS inequality (see, for example, \cite{book}) and the fact
that $\sigma(0)=1$ is the unique maximum of the continuous on $[-S,S]$
function $\sigma(t)$ imply that for some $a>1/2$, $b>0$, and all positive $%
u,\varepsilon$,
\begin{equation*}
P(E\setminus\lbrack-\varepsilon,\varepsilon];u)\leq b\exp\left(
-au^{2}\right) .
\end{equation*}
By assumptions $\mathbf{A3},\mathbf{A4}$, for the standardized process $\bar{%
X}(t)=X(t)/\sigma(t),$ $t\in\lbrack-S,S]$, for any small enough $%
\varepsilon>0$ (hence $\sigma(t)>0,$ $t\in\lbrack-\varepsilon ,\varepsilon]$%
), and for any $s,t\in\lbrack-\varepsilon,\varepsilon]$, the
following relation holds
\begin{equation*}
\mathbf{E}(\bar{X}(s)-\bar{X}(t))^{2}=2(1-\rho(s,t))\leq c_{0}|t-s|^{\gamma}
\end{equation*}
where $c_{0},\gamma$ are some positive values. Applying Theorem 8.1, \cite%
{book}, to $\bar{X}$ and definition of $B_{u}$, we obtain, that
\begin{align*}
P(\sup_{t\in E\cap\lbrack-\varepsilon,\varepsilon]\setminus B_{u}}\bar {X}%
(t)\sigma(t) & >u)\leq P(\sup_{t\in\lbrack-\varepsilon,\varepsilon]}\bar{X}%
(t)>u/\sqrt{1-u^{-2}\log^{A}u}) \\
& \leq C_{1}u^{c_{1}}\exp\left( -\frac{u^{2}}{2-2u^{-2}\log^{A}u}\right) \\
& \leq C_{2}u^{c_{1}}\exp\left( -c_{2}\log^{A}u\right) \exp\left( -\frac{%
u^{2}}{2}\right) ,
\end{align*}
for some positive $c_{i},C_{i},i=1,2.$ Since $0\in E$ by assumption, $%
P(E;u)\geq P(X(0)>u)=\Psi(u)$ for any $u>0$. Hence, the claim follows for
any $B\in(1,A)$.

\subsection{Stationary-like case}

Consider first the stationary-like case which generalize the case $%
\beta>\alpha$ in notation of \cite{PP}, \cite{book}, \cite{lectures}. Denote
for any $t\in [-S,S]$,
\begin{equation*}
f(t)=\frac{1}{2}(1-\sigma^{2}(t)),
\end{equation*}
and introduce the monotone rearrangements $f_{+}(t)$ and $f_{-}(t)$ for $%
f(t),$ $t\in\lbrack0,S]$ and $f(t),$ $t\in\lbrack-S,0],$ respectively, which
are defined as the generalized inverses
\begin{equation*}
f_{\pm}=F_{\pm}^{\leftarrow},
\end{equation*}
where
\begin{equation*}
F_{+}(x)=\mathrm{mes}\{x:f(t)\leq x,t\in\lbrack0,S]\},x\in\lbrack0,1],
\end{equation*}
and
\begin{equation*}
F_{-}(x)=\mathrm{mes}\{x:f(t)\leq x,t\in\lbrack-S,0]\},x\in\lbrack0,1],
\end{equation*}
are the distribution functions for the corresponding occupation measures,
see, for example \cite{GermanHorowitz}.

An important property of monotone rearrangements is that for any monotone
function $\phi$ we have
\begin{equation}
\int_{0}^{S} \phi(f(t)) dt = \int_{0}^{S} \phi(f_{+}(t)) dt,  \label{eq:mo}
\end{equation}
and similar equality holds for $f_{-}$.

\begin{remark}
\label{monotone} If $\sigma(t)$ is locally monotone at zero from both sides,
then for some $\varepsilon>0,$ $f_{-}(t)=f(t),$ $t\in\lbrack-\varepsilon,0],$
and $f_{+}(t)=f(t),$ $t\in\lbrack0,\varepsilon],$ i.e. for $x\in\lbrack0,1]$%
,  $F_{\pm}(x)=f^{\leftarrow}(x).$
\end{remark}

Lemma \ref{extracting} implies that the distribution functions $F_{+}(x)$
and $F_{-}(x),$ $x\in\mathbb{R}_{+}$, may be defined outside $%
[0,u^{-2}\log^{A}u],$ see (\ref{B_u}), in arbitrary way, and the asymptotic
behavior of $P([-S,S];u)$ will remain the same.

Let us introduce the Laplace transforms,%
\begin{equation}
L_{f_{+}}(\lambda):=\int_{0}^{1}e^{-\lambda x}dF_{+}(x)\ \ \text{and }%
L_{f-}(\lambda):=\int_{0}^{1}e^{-\lambda x}dF_{-}(x) , \quad\lambda>0.
\label{Lfg}
\end{equation}

\begin{theorem}
\label{stationary like d=1} Under the conditions of Lemma
\ref{extracting} together with the equality
$h_{1}(t)=0,t\in\lbrack-S,S]$, we have,
\begin{equation}
P([0,S],u)=H_{\alpha}L_{f_{+}}(u^{2})\ q^{-1}(u)\Psi (u)(1+o(1)),\ \
\label{case S}
\end{equation}
and
\begin{equation}
P([-S,S],u)=H_{\alpha}\left( L_{f_{+}}(u^{2})+L_{f_{-}}(u^{2})\right)
q^{-1}(u)\Psi(u)(1+o(1))\ \   \label{case S2}
\end{equation}
as $u\rightarrow\infty$.
\end{theorem}

\textbf{Proof}: First we consider a simplified model for $X$ and
then use Slepian inequality to derive the result for general $X$,
this is a standard approach, see \cite{book}. Let
$X_{0}(t),t\in\lbrack-S,S]$, be a centered stationary Gaussian
process satisfying conditions of Theorem \ref{Pickands}. Suppose
for a while that
\begin{equation*}
X(t)=X_{0}(t)\sigma(t),\ \ t\in\lbrack-S,S],
\end{equation*}
so that $X(t)$ satisfies the assumptions $\mathbf{A1-A4.}$ Recall that we
consider the case
\begin{equation}
\lim_{t\rightarrow0}\frac{1-\sigma^{2}(t)}{1-\rho(t)}=0,  \label{case S1}
\end{equation}
regardless of the sign of $t.$ Let us denote
\begin{equation}
T_{+}(u)=\sup\{t:t\in B_{u}\},\ \ T_{-}(u)=\inf\{t:t\in B_{u}\}.
\label{T+T-}
\end{equation}
Obviously, $T_{+}(u)>0$, $T_{-}(u)<0,$ and both of them tend to zero as $%
u\rightarrow\infty.$ In the case of locally both sides monotone $\sigma(t)$,
$T_{-}$ and $T_{+}$ are negative and positive solutions of the equation $%
1-\sigma^{2}(t)=u^{-2}\log^{A}u$, respectively, provided $u$ is sufficiently
large. Now denote
\begin{equation*}
\kappa_{-}(u):=\sqrt{q(u)T_{-}(u)},\ \ \text{and\ \ }\kappa_{+}(u):=\sqrt {%
q(u)T_{+}(u)}.\text{ }
\end{equation*}
By \eqref{q(u)} and \eqref{case S1}, and the definition of $B_{u}$,
\begin{equation*}
\lim_{u\rightarrow\infty}\frac{T_{\pm}(u)}{q(u)}=\infty,
\end{equation*}
implying
\begin{equation*}
\lim_{u\rightarrow\infty}\frac{\kappa_{\pm}(u)}{q(u)}=\lim_{u\rightarrow
\infty}\frac{T_{+}(u)}{\kappa_{+}(u)}=\lim_{u\rightarrow\infty}\frac{T_{-}(u)%
}{\kappa_{-}(u)}=\infty.
\end{equation*}

The functions $\kappa_{\pm}(u)$ satisfy the conditions of Theorem \ref%
{Pickands}. Hence, for $X_{0}(t),t\in\lbrack-S,S]$, Theorem \ref{Pickands}
implies
\begin{equation}
P\Bigl(\max_{t\in\lbrack0,\kappa_{+}(u)]}X_{0}(t)>u\Bigr)=(1+\gamma
_{+}(u))\kappa_{+}(u)H_{\alpha}q^{-1}(u)\Psi(u),  \label{P(A0)}
\end{equation}
and
\begin{equation}
P\Bigl(\max_{t\in\lbrack-\kappa_{-}(u),0]}X_{0}(t)>u\Bigr)=(1+\gamma
_{-}(u))\kappa_{-}(u)H_{\alpha}q^{-1}(u)\Psi(u),  \label{P(A0-)}
\end{equation}
where $\gamma_{\pm}(u)\rightarrow0$ as $u\rightarrow\infty.$ Denote
\begin{equation*}
\Delta_{k}(u)=k\kappa(u)+[0,\kappa(u)\mathbf{],}\quad k\in\mathbb{Z},u>0,
\end{equation*}
where we write $\kappa(u)$ instead of $\kappa_{+}(u)$ and $\kappa_{-}(u)$
and the corresponding sign depends on that side from zero (right or left) to
which $\Delta_{k}$ belongs. For all $k$ with $\Delta_{k}(u)\cap
B_{u}\neq\varnothing$, introduce the events
\begin{equation*}
A_{k}(u)=\Bigl\{\max_{t\in\Delta_{k}(u)}X_{0}(t)>u_{k}\Bigr\},\ \text{where }%
u_{k}=u/\sigma_{k},\ \sigma_{k}=\max_{t\in\Delta_{k}(u)}\sigma(t),
\end{equation*}
and
\begin{equation*}
A_{k}^{\prime}(u)=\Bigl\{\max_{t\in\Delta_{k}(u)}X_{0}(t)>u_{k}^{\prime }%
\Bigr\},\ \text{where }u_{k}^{\prime}=u/\sigma_{k}^{\prime},\ \sigma
_{k}^{\prime}=\min_{t\in\Delta_{k}(u)}\sigma(t).
\end{equation*}
Since
\begin{equation}
u\leq u_{k},u_{k}^{\prime}\leq u+u^{-1}\log^{A}u,\ k\in K_{u}:=\{k:\Delta
_{k}(u)\cap B_{u}\neq\varnothing\},  \label{K_u}
\end{equation}
and all the intervals $\Delta_{k}(u)$ have length $\kappa(u),$ $\kappa(u)$
also satisfies the conditions of Theorem \ref{Pickands} with $u_{k}$ instead
of $u.$ Therefore,
\begin{equation*}
{P}(A_{k})=(1+\gamma(u_{k}))H_{\alpha}\kappa(u)q^{-1}(u_{k})\Psi(u_{k}),
\end{equation*}
and
\begin{equation*}
{P}(A_{k}^{\prime})=(1+\gamma(u_{k}^{\prime}))H_{\alpha}%
\kappa(u)q^{-1}(u_{k}^{\prime})\Psi(u_{k}^{\prime}).
\end{equation*}
By definition (\ref{K_u}) of $K_{u},$ since $q$ is non-increasing, there
exists a positive non-increasing $\delta_{1}(u)$ tending to zero as $%
u\rightarrow\infty,$ such that
\begin{equation*}
1-\delta_{1}(u)\leq\frac{\min(q(u_{k}),q(u_{k}^{\prime}))}{q(u)} \leq \frac{%
\max(q(u_{k}),q(u_{k}^{\prime}))}{q(u)} \leq1,\ \ \forall k\in K_{u}.
\end{equation*}
Further, since $\delta_{2}(u)=\sup_{v\geq u}|\gamma(v)|\rightarrow0$ as $%
u\rightarrow\infty$ and $u_{k},u_{k}^{\prime}\geq u,$ we obtain,
that
\begin{equation}
{P}(A_{k}),{P}(A_{k}^{\prime})\lesseqgtr(1\pm\delta_{1}(u))(1\pm\delta
_{2}(u))H_{\alpha}\kappa(u)q^{-1}(u_{k})\Psi(u_{k}).  \label{P(A_ku)}
\end{equation}
Due to Bonferroni inequalities,%
\begin{equation}
{P(B}_{u};u)\leq\sum_{k:\Delta_{k}(u)\cap
B_{u}\neq\varnothing}{P}(A_{k}(u)), \label{DSabove}
\end{equation}
and%
\begin{equation}
{P(B}_{u};u)\geq\sum_{k::\Delta_{k}(u)\subset B_{u}}{P}(A_{k}^{\prime
}(u))-\sum_{k,l:k:\Delta_{k}(u)\cap B_{u}\neq\varnothing,\Delta_{l}(u)\cap
B_{u}\neq\varnothing,k\neq l}{P}(A_{k}(u)A_{l}(u)).  \label{DSbelow}
\end{equation}
We do not assume symmetry property of $\sigma(t),$ therefore consider
separately the sums in (\ref{DSabove}) with positive and negative $k.$ Let
us rewrite
\begin{equation}
u_{k}^{2}=u^{2}+u^{2}(1-\sigma_{k}^{2})+\frac{u^{2}(1-\sigma_{k}^{2})^{2}}{%
\sigma_{k}^{2}}.  \label{u_k}
\end{equation}
The inequality
\begin{equation*}
\frac{u^{2}(1-\sigma_{k}^{2})^{2}}{2\sigma_{k}^{2}}\geq\frac{%
u^{2}u^{-4}\log^{2A}u}{2},
\end{equation*}
implies that for some positive $\delta_{3}(u)$ with $\delta_{3}(u)%
\rightarrow0,$ $u\rightarrow\infty,$ uniformly in $k\in K_{u},$%
\begin{equation*}
\sigma_{k}\exp\left( -\frac{u^{2}(1-\sigma_{k}^{2})^{2}}{2\sigma_{k}^{2}}%
\right) \lesseqgtr1\pm\delta_{3}(u).
\end{equation*}
Hence, by (\ref{P(A_ku)}),
\begin{equation*}
\sum_{k:\Delta_{k}(u)\cap B_{u}\neq\varnothing}{P}(A_{k}(u)) \lesseqgtr
(1\pm\delta_{4}(u)) H_{\alpha} \kappa(u) q(u)^{-1}\Psi(u)\sum_{k\geq 0:\
\Delta_{k}\cap B_{u}\neq\varnothing}e^{-u^{2}(1-\sigma_{k}^{2})/2},
\end{equation*}
where
\begin{equation*}
1\pm\delta_{4}(u)=(1\pm\delta_{1}(u))(1\pm\delta_{2}(u))(1\pm\delta_{3}(u)).
\end{equation*}
Define
\begin{equation*}
\Sigma(u)=\sum_{k\geq0:1-\sigma_{k}^{2}\leq
u^{-2}\log^{A}u}e^{-u^{2}(1-\sigma_{k}^{2})/2}.
\end{equation*}
Hence, $\kappa(u)\Sigma(u)$ is an integral sum for the integral%
\begin{equation*}
I_f(u)=\int_{0\leq f(t)\leq2u^{-2}\log^{A}u}e^{-u^{2}f(t)}dt.
\end{equation*}
By \eqref{eq:mo},
\begin{equation*}
I_f(u)=I_{f+}(u)=\int_{0\leq
f_{+}(t)\leq2u^{-2}\log^{A}u}e^{-u^{2}f_{+}(t)}dt.
\end{equation*}
Change the order of summands in $\Sigma(u)$ by ordering $u_{k}$'s (or
equivalently $\sigma_{k}$'s) gives that $\kappa(u)\Sigma(u)$ is an integral
sum for $I_{+}(u)$ as well. Since $f_{+}(t)$ is monotone in the integration
domain for sufficiently large $u$, we have,
\begin{equation*}
\kappa(u)(\Sigma(u)-1)\leq I_f(u)<\kappa(u)(\Sigma(u)+1).
\end{equation*}
The limit $\lim_{u \to\infty} \kappa(u)/q(u)= \infty$ implies that
\begin{equation*}
\lim_{u \to\infty} \Sigma(u)=\infty,
\end{equation*}
and consequently,
\begin{equation*}
I_f(u)=(1+o(1)) \kappa(u)\Sigma(u)
\end{equation*}
as $u\to\infty$. Recall that if $\sigma(t)$ is locally monotone at zero then
for all sufficiently large $u,$ $f_{+}(t)=f(t)$ in the integration domain.
Changing variables $x=$ $f_{+}(t)$ and denoting $\lambda=u^{2}$ implies
\begin{equation*}
I_{f+}(u)\sim\int_{0}^{2\lambda^{-1}\log^{A/2}\lambda}e^{-\lambda x}dF_{+}(x)
\end{equation*}
as $\lambda\rightarrow\infty.$ Consideration of the parts over the negative
values of $k$ in the single sums is the same. Thus, for the single sums in (%
\ref{DSabove}, \ref{DSbelow}), we get
\begin{align*}
& \sum_{k::\Delta_{k}(u)\cap B_{u}\neq\varnothing}{P}(A_{k}(u))=(1+o(1))\sum
_{k:\Delta_{k}(u)\subset B_{u}}{P}(A_{k}^{\prime}(u)) \\
& =(1+o(1))\mathcal{H}_{\alpha}\left( I_{f-}(u)+I_{f+}(u)\right)
q^{-1}(u)\Psi(u)
\end{align*}
as $u\rightarrow\infty.$ The estimation of the double sum is quite similar
to that in \cite{book}, \cite{lectures}.

Further, for {$\lambda =u^{2}$} and any $A>1$,
\begin{equation*}
\int_{2\lambda ^{-1}\log ^{A/2}\lambda }^{1}e^{-\lambda x}dF(x)\leq \exp
\left( -2\log ^{A/2}u^{2}\right) .
\end{equation*}%
Since $P([0,S],u)\geq \Psi (u)$ the above asymptotics of $P([0,S],u)$ and
the fact that $A$ can be taken arbitrary large, implies that
\begin{equation}
\lim_{u\rightarrow \infty }\frac{I_{f\pm }(u)}{L_{f\pm }(u^{2})}=1.
\label{eq:I}
\end{equation}

Hence in view of already mentioned standard passage from the particular $%
X(t)=X_{0}(t)\sigma(t)$ to the general Gaussian process (by applying Slepian
inequality), the proof follows easily.

\subsubsection{Remarks on representations and properties of $L_{f_{\pm}}(%
\protect\lambda).$}

Our assumptions give that $\sigma $ is continuous and attains its unique
maximum on $[-S,S]$ at 0. Consequently, for such $\sigma $ as shown above,
the relation \eqref{eq:I} holds. First, consider an important case of
regularly varying $F_{\pm }(x)$ at zero, i.e., for some slowly varying at
zero functions $\ell _{\pm }(t)$,%
\begin{equation}
F_{\pm }(x)\sim \ell _{\pm }(x)x^{a_{\pm }}\ \ \text{as }x\rightarrow 0,
\label{RV}
\end{equation}%
where $a_{\pm }\geq 0$. By \cite{feller}, Theorems XIII.2 and XIII.3, %
\eqref{RV} is equivalent to
\begin{equation*}
L_{f_{\pm }}(\lambda )\sim \Gamma (1+a_{\pm })\lambda ^{-a_{\pm }}\ell _{\pm
}(1/\lambda )\ \ \text{as }\lambda \rightarrow \infty .
\end{equation*}

Using (\ref{q(u)}), we obtain the following statement.

\begin{proposition}
\label{stationary like PP}Suppose that the above assumptions hold and $%
h_{1}(t)\equiv0.$ If (\ref{RV}) holds with $\ell_{\pm}, a_{\pm}$ as above,
then $a_{\pm}\le 1/\alpha$ and
\begin{equation}
P([0,S],u)=H_{\alpha}\Gamma(1+a_{+})u^{ \frac{2}{\alpha}-2a_{+}}%
\ell_{+}(u^{2})(\ell^{\#}(u^{-2}))^{\frac{1}{\alpha}}\Psi(u)(1+o(1)),
\label{RV1}
\end{equation}

\begin{eqnarray}
&P([-S,S],u)=H_{\alpha}\left( \Gamma(1+a_{+})u^{\frac{2}{\alpha }%
-2a_{+}}\ell_{+}(u^{2})+ \Gamma(1+a_{-})u^{\frac{2}{\alpha}-2a_{-}}\ell
_{-}(u^{2})\right)  \notag \\
&\times (\ell^{\#}(u^{-2}))^{\frac{1}{\alpha}}\Psi (u)(1+o(1)),  \label{RV2}
\end{eqnarray}
as $u\rightarrow\infty$.
\end{proposition}

\begin{remark}
\label{BGT} i) If $f_{\pm}(x)= x^{\beta_{\pm}}\ell_{\pm}(x) $ are regularly
varying at zero with positive indexes $\beta_{\pm}$ ( $\ell_{\pm}(x)$ are
slowly varying), then by Theorem 1.5.3, \cite{BGT} there exists an
asymptotically monotone equivalent to $f_{\pm}$, say, $f_{*,\pm}(t)$. Hence,
in view of the above argument, one can take $F_{\pm}(x)=f_{*,\pm}^{%
\leftarrow }(x).$ Further, by the same argument as before in (\ref{q(u)}),
we have
\begin{equation*}
F_{\pm}(x)\sim x^{1/\beta_{\pm}}\ell^{\#}_{\pm}(x),\ \ x\rightarrow0.
\end{equation*}
ii) The case when $f$ is regularly varying at $0$ has been recently
investigated in \cite{HDL}, where the authors established \eqref{RV1}. In
fact, $1-\sigma $ is assumed to be symmetric around 0 therein, the
non-symmetric case can be established with no additional efforts. Note
further that the case $\ell_{\pm }(x)=1,x\in[-S,S]$ was considered in \cite%
{PP}, in this case $L_{f_{\pm}}(\lambda)\sim\Gamma(1+1/\beta)\lambda^{-1/%
\beta}$.
\end{remark}

Now consider some other representations of $L_{f_{\pm}}(\lambda).$
Integration by parts and choice of sufficiently large $A$ imply that in the
stationary-like case for any $\lambda >0$,
\begin{eqnarray*}
L_{f\pm}(\lambda)&=e^{-2\log^{A/2}\lambda}+\lambda\int_{0}^{2\lambda^{-1}%
\log^{A/2}\lambda}e^{-\lambda x}F_{+}(x)dx \\
&=: e^{-2\log^{A/2}\lambda}+ J_{f_{\pm}}(\lambda) \sim J_{f_{\pm}}(\lambda)
\end{eqnarray*}
as $\lambda\to\infty$. Moreover, if $\sigma(t)$ is locally monotone at $0$,
\begin{equation}
J_{f_{\pm}}(\lambda)\sim\int_{\exp(-2\log^{A/2}\lambda)}^{1}f^{\leftarrow
}(-\lambda^{-1}\log v)dv, \quad\lambda\to\infty.  \label{Afm}
\end{equation}

\subsection{Talagrand case}

It immediately follows from Lemmas \ref{PP-lemma} and \ref{extracting} that
if $h_{1}(t)=\infty$ for all $t\neq0,$ then
\begin{equation}
P([0,S],u)=P([-S,S],u)=\Psi(u)(1+o(1)),\ \ \   \label{case T}
\end{equation}
as $u\rightarrow\infty.$

\subsection{The transition case}

We already know that in this case $1-\sigma^{2}(t)=Ct^{\alpha}\ell_{1}(t)$
with $\ell(t)/\ell_{1}(t)\rightarrow1$ as $t\rightarrow0,$ where $%
1-r(s,t)\sim|t-s|^{a}\ell(t-s),$ $s,t\rightarrow0.$ In this case the
exceeding probability asymptotic evaluation is very similar to the
corresponding evaluations in \cite{PP}, \cite{book}, \cite{lectures}. The
only difference is that the case $\ell(t)=\ell_{1}(t)=1$ is considered in
these papers. As shown in \cite{HDL}, the slowly varying function does not
play any role in the asymptotics. Consequently, in this case we obtain
\begin{equation}
P([0,S],u)=(1+o(1))P_{\alpha}^{+}\Psi(u),  \label{case P}
\end{equation}%
\begin{equation}
P([-S,S],u)=(1+o(1))P_{\alpha}\Psi(u),  \label{case PP}
\end{equation}
as $u\rightarrow\infty,$ where $P_{\alpha}^{+}=\lim_{T\rightarrow\infty
}P_{\alpha}(T)\in(0,\infty),\ $and $\mathbf{\ }P_{\alpha}^{+}=\lim
_{T\rightarrow\infty}P_{\alpha}(T)\in(0,\infty).$

Note that in contrast with the stationary like case, in the Talagrand and
transition cases the double side probabilities are not asymptotically equal
to the sum of one side ones.

\subsection{Main result}

Now we combine all the obtained results concerning the asymptotic behavior
in the non-stationary case and formulate our main result. We say that we have

\begin{itemize}
\item \emph{S-S case,} the stationary like case (considered in Proposition %
\ref{stationary like d=1});

\item \emph{S-T case,} when $h_{1}(t)=0,$ $t\leq0,$ and $h_{1}(t)=\infty,$ $%
t>0;$

\item \emph{P-S case,} when $h_{1}(t)\in(0,\infty),$ $t\leq0,$ and $%
h_{1}(t)=0,$ $t>0;$

\item so on, similarly for the remaining $6$ cases.
\end{itemize}

\begin{theorem}
\label{main} If $X(t),$ $t\in\lbrack-S,S]$ is a Gaussian zero mean process
satisfying conditions $\mathbf{A1-A4}$, then:

\begin{itemize}
\item In \emph{S-S case}, (\ref{case S2}) is valid;

\item In other four cases concerning \emph{S} the asymptotic behavior is
equal to the right hand side of (\ref{case S});

\item In \emph{T-T case}, (\ref{case T}) is valid;

\item In \emph{P-P case,} (\ref{case PP}) is valid;

\item In \emph{T-P, P-T cases} the asymptotic behavior is equal to the right
hand side of (\ref{case P}).
\end{itemize}
\end{theorem}

\begin{remark}
{i) If a set $E\subset [-S,S]$ is a closure of an open bounded set
containing the unique point of the variance maximum then $P(E,u) \sim
P([-S,S],u)$ as $u\to\infty$. In particular, all the asymptotic results
above do not depend on $S>0$. }

ii) In the case when the maximum point of $\sigma^{2}$ is a boundary point,
the corresponding one side relations hold.
\end{remark}

\section{Examples}

Below we present two illustrating examples of \emph{S-S case}. Exotic cases
when $\sigma(t)$ is not locally monotone can be dealt similarly by
calculating first the monotone rearrangement of $f(t)=1-\sigma^{2}(t)$.

\textbf{Example 1.} The case of very gentle sharp maximum. Let for some $%
\varepsilon,\beta>0,$ and positive $t,$
\begin{equation*}
f(t)=1-\sigma^{2}(t)=e^{-t^{-\beta}},\ \ t\in(0,\varepsilon].
\end{equation*}
We have that $\sigma(t)$ is symmetric and locally monotone on both sides, so
we write simply $F_{\pm}=F.$ For $x>0,$
\begin{equation*}
F(x)=f^{\leftarrow}(x)=\log^{-1/\beta}(1/x),
\end{equation*}
that is ${a_{\pm}}=0,$ $\ell(x)=\log^{-1/\beta}(1/x).$ Hence by {Proposition}
\ref{stationary like PP},%
\begin{equation*}
P([-S,S],u)=2^{1-1/\beta}H_{\alpha}\log^{-1/\beta}(u) \frac{\Psi (u)}{{q(u)}}%
(1+o(1))
\end{equation*}
as $u\rightarrow\infty.$\newline
The asymptotic behaviors of $P([-S,S],u)$ for this and similar "gentle"
cases  has not been known in the literature so far.

\textbf{Example 2.} Consider the case when $1-\sigma^{2}(t)$ is close to $%
1-r(s,t),$ $s,t\rightarrow0.$ Take $1-\sigma^{2}(t)=t^{\alpha}\log(1/t),$
with, of course, $\log(1/t)/\ell(t)\rightarrow0,$ as $t\rightarrow0,$ where,
combining $\mathbf{A3}$ and $\mathbf{A4,}$ $1-r(s,t)\sim|t-s|^{a}\ell(t-s),$
$s,t\rightarrow0.$ By Corollary 2.3.4, \cite{BGT}, see also Remark \ref{BGT}%
,
\begin{equation*}
(\log(1/t))^{\#}\sim\frac{1}{\log(1/t)},
\end{equation*}
so that
\begin{equation*}
F_{\pm}(x)\sim\frac{x^{1/\alpha}}{\log(1/x)},\ \ x\rightarrow0,
\end{equation*}
and
\begin{equation*}
L_{f_{\pm}}(u^{2})\sim\frac{\Gamma(1+1/\alpha)u^{-2/\alpha}}{2\log u}\ \
\text{as }u\rightarrow\infty.
\end{equation*}
Thus$,$%
\begin{equation*}
P([-S,S],u)=H_{\alpha}\frac{\ell^{\#}(u^{-2})^{1/\alpha}}{\log u}%
\Psi(u)(1+o(1))
\end{equation*}
as $u\rightarrow\infty,$ see (\ref{q(u)}).

\end{document}